\documentclass{amsart}
\usepackage{amssymb,amscd,amsthm}
\usepackage{latexsym}
\date{\today}

\newcommand{\Z}{{\mathbb Z}}
\newcommand{\R}{{\mathbb R}}
\newcommand{\C}{{\mathbb C}}

\newcommand{\T}{{\mathbb T}}
\newcommand{\Q}{{\mathbb Q}}

\newtheorem{theorem}{Theorem}

\newtheorem{lemma}{Lemma}

\newtheorem{coro}{Corollary}

\newtheorem{definition}{Definition}
\begin{document}
\title[The Repetition Property for Polynomials and Skew-Shifts]{The Repetition
Property for Sequences on Tori Generated by Polynomials or
Skew-Shifts}

\author[M.\ Boshernitzan]{Michael Boshernitzan}

\address{Department of Mathematics, Rice University, Houston, TX~77005, USA}

\email{michael@rice.edu}

\author[D.\ Damanik]{David Damanik}

\address{Department of Mathematics, Rice University, Houston, TX~77005, USA}

\email{damanik@rice.edu}

\thanks{D.\ D.\ was supported in part by NSF grant
DMS--0653720.}

\begin{abstract}
The repetition property of a sequence in a metric space, a notion
introduced by us in an earlier paper, is of importance in the
spectral analysis of ergodic Schr\"odinger operators. It may be
used to exclude eigenvalues for such operators. In this paper we
study the question of when a sequence on a torus that is generated
by a polynomial or a skew-shift has the repetition property. This
provides classes of ergodic Schr\"odinger operators with
potentials generated by skew-shifts on tori that have, contrary to
earlier belief, no eigenvalues.
\end{abstract}

\maketitle

\section{Introduction}

This short paper contains results that are complementary to those
obtained by us in a recent paper \cite{bd1}. There, the following
notion was introduced. (Here and in what follows, we write $\Z_+ =
\{ 1 , 2 , 3 , \ldots \}$.)

\begin{definition}
A sequence $\{\omega_n\}_{n \ge 0}$ in a metric space $\Omega$ has
the repetition property if for every $\varepsilon > 0$ and $r \in
\Z_+$, there exists $q \in \Z_+$ such that $\mathrm{dist}(\omega_n
, \omega_{n+q}) < \varepsilon$ for $n = 0, 1, 2 , \ldots , rq $.
\end{definition}

An alternative way to state this definition that involves fewer
parameters and is sometimes more convenient to use is the
following: A sequence $\{\omega_n\}_{n \ge 0}$ in a metric space
$\Omega$ has the repetition property if and only if for every $r
\in \Z_+$, there exists $q \in \Z_+$ such that
$\mathrm{dist}(\omega_n , \omega_{n+q}) < 1/r$ for $n = 0, 1, 2 ,
\ldots , rq $.

Sequences of particular interest are given by forward orbits
corresponding to some map $T : \Omega \to \Omega$, that is,
$\mathcal{O}_+(\omega) = \{ T^n \omega \}_{n \ge 0}$. For a given
map $T$, the set $PRP(T)$ consists of the points $\omega \in
\Omega$ for which $\mathcal{O}_+(\omega)$ has the repetition
property. Following \cite{bd1}, we say that $T$ has the
topological repetition property {\rm (TRP)} if $PRP(T) \not=
\emptyset$, the metric repetition property {\rm (MRP)} relative to
some $T$-invariant measure $\mu$ if $\mu(PRP(T))
> 0$, and the global repetition property {\rm (GRP)} if
$PRP(T) = \Omega$. Obviously, we always have (GRP) $\Rightarrow$
(MRP) (relative to any non-zero $T$-invariant measure)
$\Rightarrow$ (TRP). Also, if $T$ is continuous and $(\Omega,T)$
is minimal, then (TRP) implies that $PRP(T)$ is residual; compare
Corollary~\ref{c.gdelta} below.

If $(\Omega,T)$ is a symbolic flow, the repetition property for a
forward orbit $\mathcal{O}_+(\omega)$ is related to the initial
critical exponent; see \cite{bhz} for a precise definition and
further references. Namely, $\omega \in PRP(T)$ if and only if the
initial critical exponent of $\omega$ is infinite, that is,
$\omega$ has arbitrarily high powers as prefixes. For the
prominent class of Sturmian symbolic flows, (TRP) and (MRP)
relative to the unique shift-invariant probability measure are
equivalent and they hold if and only if the so-called slope of the
Sturmian flow has unbounded partial quotients. It is also of
interest to note that (GRP) fails whenever the symbolic flow
$(\Omega,T)$ is aperiodic; compare \cite{bd1}. The repetition
property in the symbolic setting and related issues will be
studied further in \cite{bd2}.

The repetition property for forward orbits plays an important role
in the study of the eigenvalue problem for Schr\"odinger operators
in $\ell^2(\Z)$,
$$
(H_\omega \psi)(n) = \psi(n+1) + \psi(n-1) + V_\omega(n) \psi(n),
$$
with dynamically defined potentials, that is, $V_\omega (n) =
f(T^n \omega)$ with an invertible map $T : \Omega \to \Omega$ and
a continuous function $f : \Omega \to \R$. Namely, in the
particular situation where $\Omega$ is compact and $T$ is a
minimal homeomorphism, the size of the set $PRP(T)$ was related in
\cite{bd1}, for a generic continuous function $f$, to the size of
the set of $\omega$'s for which $H_\omega$ has no eigenvalues: If
$T$ has (TRP) (resp., (MRP) relative to some $T$-ergodic measure
$\mu$), then there is a residual set $\mathcal{F} \subset
C(\Omega)$ such that for each $f \in \mathcal{F}$, $H_\omega$ has
no eigenvalues for $\omega$'s from a residual (resp., full
$\mu$-measure) subset of $\Omega$.

Several examples were considered in \cite{bd1}. For most of them,
$\Omega$ is either the circle $\T$ or a $k$-dimensional torus
$\T^k$. Here, we denote $\T = \R / \Z$. The metric on $\T^k$ is
given by $\mathrm{dist}(x,y) = \sum_{j = 1}^k \langle x_j - y_j
\rangle$, where $\langle \tau \rangle = \min \{ |\hat \tau - p| :
p \in \Z \}$ with any representative $\hat \tau \in \R$ of $\tau
\in \T$. For a shift on a torus (i.e., $T\omega = \omega +
\alpha$), it is easily seen that (GRP) holds. Moreover, we showed
that almost every interval exchange transformation has (MRP)
relative to Lebesgue measure. On the other hand, almost every
interval exchange transformation is weakly mixing.\footnote{The
permutation of intervals is assumed to be irreducible and not a
rotation.} This is a deep result with significant contributions by
Katok, Stepin, Veech, Nogueira, Rudolph and others; the final
result is due to Avila and Forni \cite{af}. This provides a class
of examples that have (MRP) and are weakly mixing.

As explained above, one obtains as a consequence results on the
absence of eigenvalues for Schr\"odinger operators with potentials
generated by shifts on a torus or by an interval exchange
transformation that hold for almost every point in the space in
question. To a certain extent this confirms what was expected for
such operators based on earlier results.

For Schr\"odinger operators with potentials generated by a
skew-shift on a torus, however, it was expected that there should
be plenty of eigenvalues; see, for example, Bourgain's recent book
\cite{b5} and the discussion in \cite{bd1}. Alas, we were able to
give a complete study of (TRP), (MRP), and (GRP) for skew-shifts
on the two-dimensional torus, $(\omega_1,\omega_2) \mapsto
(\omega_1 + \alpha , \omega_2 + \omega_1)$, and proved in
particular that all three properties hold for Lebesgue almost
every $\alpha$. Hence we obtained the corresponding absence of
eigenvalue result for the associated Schr\"odinger operators which
was surprising given earlier expectations for these operators.
Possible extensions to higher-dimensional tori were only
indicated. One of the main goals of the present paper is to work
out these extensions in detail.

The class of skew-shifts we are interested in is given by maps $T
: \T^k \to \T^k$, $k \ge 2$, that have an irrational rotation in
the first component and, for $j \ge 2$, the $j$-th component is
given by a rational linear combination of $\omega_1, \ldots,
\omega_j$. For example, for $\alpha$ irrational, the following map
fits this profile:
\begin{equation}\label{skewshift}
T : \T^k \to \T^k, \; (\omega_1 , \omega_2 , \ldots , \omega_k)
\mapsto ( \omega_1 + \alpha , \omega_2 + \omega_1 , \ldots ,
\omega_k + \omega_{k-1} ).
\end{equation}
We refer the reader to the paper \cite{f} by Furstenberg, which
studies skew-shifts of this (and more general) form from an ergodic
theory point of view. It was shown there that normalized Lebesgue
measure on $\T^k$ is the unique $T$-invariant probability measure.
Thus, in the context of the map \eqref{skewshift}, (MRP) will always
be understood to be relative to this measure.

While our arguments will also apply to more general skew-shifts on
the torus, we will consider for definiteness the specific map
\eqref{skewshift}. It turns out that the size of $PRP(T)$ depends
sensitively on both the dimension $k$ of the torus and Diophantine
properties of the irrational number $\alpha$. The following result
gives almost complete information as to which of the the
properties (TRP), (MRP) (relative to Lebesgue measure), and (GRP)
hold.

\begin{theorem}\label{t.skewshift}
Suppose that $\alpha \in \T$ is irrational and the map $T : \T^k
\to \T^k$ is given by $T(\omega_1 , \omega_2 , \ldots , \omega_k)
= ( \omega_1 + \alpha , \omega_2 + \omega_1 , \ldots , \omega_k +
\omega_{k-1} )$.
\\
{\rm (a)} For every $k \ge 2$, the following are equivalent:
\begin{itemize}

\item $T$ has {\rm (TRP)},

\item $\liminf_{q \to \infty} q^{k-1} \langle q \alpha \rangle =
0$.

\end{itemize}
{\rm (b)} For $k = 2$, the following are equivalent:
\begin{itemize}

\item $T$ has {\rm (TRP)},

\item $T$ has {\rm (MRP)},

\item $T$ has {\rm (GRP)},

\item $\liminf_{q \to \infty} q \langle q \alpha \rangle = 0$.

\end{itemize}
{\rm (c)} If $k \ge 3$, then $T$ does not have {\rm (GRP)}.
\\
{\rm (d)} If $k \ge 3$, then $T$ does not have {\rm (TRP)} for
Lebesgue almost every $\alpha$.
\\
{\rm (e)} If $k = 3$, then $T$ has {\rm (MRP)} for $\alpha$'s from
a residual subset of $\T$.
\\
{\rm (f)} If $k \ge 4$, then $T$ does not have {\rm (MRP)}.
\end{theorem}

In particular, parts (a) and (e) of Theorem~\ref{t.skewshift}
yield new classes of Schr\"odinger operators generated by
skew-shifts that (surprisingly) have empty point spectrum.

Since the components of $T^n \omega$ are given by the projection
to $\T$ of polynomials in $n$ with coefficients that depend on
$\alpha$ and $\omega_1, \ldots, \omega_k$, we will first study the
repetition property for sequences on the circle $\T$ that are
generated by polynomials. This is done in Section~\ref{s.poly}.
Skew-shift orbits are then investigated in Section~\ref{s.skew},
where we prove Theorem~\ref{t.skewshift}.

\section{General Preliminary Results}

This section contains several results that hold in the general
setting and which will be useful later when we specialize to the
circle or the $k$-dimensional torus.

\begin{lemma}\label{l.gdelta}
Assume that $X$ is a topological space, $\Omega$ is a metric
space, and $f_n : X \to \Omega$, $n \in \Z_+$ are continuous. Then
the set
$$
RPX = \{ x \in X : \text{ the sequence } \{ f_n(x) \}_{n \ge 0}
\text{ has the repetition property} \, \}
$$
is a $G_\delta$ set in $X$.
\end{lemma}

\begin{proof}
Denote for $r,q \in \Z_+$ and $\varepsilon > 0$,
$$
H(r,q,\varepsilon) = \{ x \in X : \max_{0 \le n \le rq }
\mathrm{dist} ( f_n(x), f_{n+q}(x) ) < \varepsilon \}.
$$
Since each set $H(r,q,\varepsilon)$ is open in $X$ and
$$
RPX = \bigcap_{r \ge 1} \bigcap_{k \ge 1} \bigcup_{q \ge k}
H(r,q,\tfrac{1}{r}),
$$
it follows that $RPX$ is a $G_\delta$ set in $X$.
\end{proof}

\begin{coro}\label{c.gdelta}
If $\Omega$ is a metric space and $T : \Omega \to \Omega$ is
continuous, then $PRP(T)$ is a $G_\delta$ subset of $\Omega$.
\end{coro}

\begin{proof}
Let $X = \Omega$, $f_n(\omega) = T^n \omega$ for $\omega \in
\Omega$ and $n \ge 0$, and apply the previous lemma.
\end{proof}

Next we recall the topological lemma from Boshernitzan's appendix
to Cheung's paper \cite{bc}. It will play a crucial role in the
the proof of Theorem~\ref{t.border} below, which in turn is an
essential ingredient to the proof of
Theorem~\ref{t.skewshift}.(e).

\begin{lemma}\label{l.chebosh}
Assume that $L$ is a $G_\delta$ subset of a $\sigma$-compact
metric space $K$, $P$ is a Polish space, $H$ is an $F_\sigma$
subset of $W = P \times K$. Denote, for $p \in P$,
$$
K(p) = \{ k \in K : (p,k) \in H \}
$$
and
$$
P_0 = \{ p \in P : K(p) \subset L \}.
$$
If $P_0$ is dense in $P$, then $P_0$ is a residual subset of $P$.
\end{lemma}

\begin{proof}
See \cite[Lemma~A.1]{bc}.
\end{proof}

The following notion will prove to be useful in our study of the
repetition property for sequences generated by polynomials and
skew-shifts.

\begin{definition}
A family of sequences $\{ \omega_n^{(\gamma)} \}_{n \ge 0}$ in
metric spaces $\Omega^{(\gamma)}$, $\gamma \in \Gamma$, has the
joint repetition property if each of them has the repetition
property and for each finite subfamily, $q = q(\varepsilon,r)$ can
be chosen uniformly for all sequences in the finite subfamily.
\end{definition}

We now state three lemmas about the joint repetition property. The
proofs are so simple that we omit them.

\begin{lemma}\label{l.joint1}
Suppose $\Omega,\Omega'$ are metric spaces, $\{\omega_n\}_{n \ge
0}$ is a sequence in $\Omega$, and $\{\omega_n'\}_{n \ge 0}$ is a
sequence in $\Omega'$. Then, the sequences $\{\omega_n\}_{n \ge
0}$ and $\{\omega_n'\}_{n \ge 0}$ have the joint repetition
property if and only if the sequence $\{(\omega_n,\omega_n')\}_{n
\ge 0}$ in $\Omega \times \Omega'$ has the repetition property.
\end{lemma}

\begin{lemma}\label{l.joint2}
Suppose $\Omega,\Omega'$ are metric spaces, $g : \Omega \to
\Omega'$ is uniformly continuous, $\{\omega_n\}_{n \ge 0}$ is a
sequence in $\Omega$ with the repetition property, and the
sequence $\{\omega_n'\}_{n \ge 0}$ in $\Omega'$ is given by
$\omega_n' = g(\omega_n)$ for $n \ge 0$. Then, the sequences
$\{\omega_n\}_{n \ge 0}$ and $\{\omega_n'\}_{n \ge 0}$ have the
joint repetition property.
\end{lemma}

\begin{lemma}\label{l.joint3}
Suppose $\Omega,\Omega'$ are metric spaces and $\{\omega_n\}_{n
\ge 0}$ is a sequence in $\Omega$ with the repetition property.
Then, the set of sequences $\{\omega_n'\}_{n \ge 0}$ in $\Omega'$
so that $\{\omega_n\}_{n \ge 0}$ and $\{\omega_n'\}_{n \ge 0}$
have the joint repetition property is closed with respect to
uniform convergence.
\end{lemma}

\section{Sequences on the Circle}\label{s.poly}

In this section we study sequences on the circle, that is, we let
$\Omega = \T$. We will be especially interested in sequences
generated by polynomials.

\subsection{Preservation Properties}

In this subsection we exhibit a number of operations performed on a
given sequence that preserve the repetition property.

We begin with operations that obviously preserve the repetition
property.

\begin{lemma}\label{l.basic}
Suppose $\{\omega_n\}_{n \ge 0}$ is a sequence in $\T$ that has the
repetition property. Then, the following sequences have the
repetition property as well:
\\
{\rm (a)} $\{\omega_{n+l}\}_{n \ge 0}$ for every $l \in \Z_+$,
\\
{\rm (b)} $\{\omega_{nl}\}_{n \ge 0}$ for every $l \in \Z_+$,
\\
{\rm (c)} $\{l \omega_n\}_{n \ge 0}$ for every $l \in \Z$.
\\
In fact, the sequences listed above have the joint repetition
property.
\end{lemma}

\begin{proof}
This is readily verified.
\end{proof}

Note that parts (a) and (b) of this lemma are not specific to the
circle and hold in a general metric space.

The next lemma addresses the following question: Which sequences can
be added to any given sequence without destroying the repetition
property? This is related to having the joint repetition property
with any given sequence that has the repetition property. Here we
only treat those sequences we need in the sequel. In the appendix we
investigate this issue in more depth.

\begin{lemma}\label{l.addition}
Suppose $\{\omega_n\}_{n \ge 0}$ is a sequence in $\T$ that has the
repetition property. Then, for every $\alpha \in \T$,
$\{\omega_n\}_{n \ge 0}$ and $\{ \alpha n \}_{n \ge 0}$ have the
joint repetition property. Consequently, for every $\alpha, \beta
\in \T$, $\{\omega_n\}_{n \ge 0}$ and $\{\omega_n + \alpha n +
\beta\}_{n \ge 0}$ have the joint repetition property.
\end{lemma}

\begin{proof}
Since $\{\omega_n\}_{n \ge 0}$ has the repetition property, there
exists, for every small $\varepsilon > 0$ and $r \in \Z_+$, an
integer $\tilde q \ge 1$ such that
$$
\langle \omega_{n + \tilde q} - \omega_n \rangle < \varepsilon^2
\quad \text{ for } n = 0, 1, 2, \ldots,
\left\lfloor\frac{r+2}{\varepsilon} \right\rfloor \tilde q.
$$
Select an integer $1 \le d \le \frac{1}{\varepsilon}$ such that
$\langle d \tilde q \alpha \rangle < \varepsilon$. Set $q = d \tilde
q$. Then we have
$$
\langle \omega_{n + q} - \omega_n \rangle < d \varepsilon^2 \le
\varepsilon \quad  \text{ for }  n = 0, 1, 2, \ldots, r q
$$
and
$$
\langle \alpha (n + q) - \alpha n \rangle = \langle q \alpha \rangle
= \langle d \tilde q \alpha \rangle < \varepsilon \quad \text{ for
every } n \in \Z.
$$
It follows that$\{\omega_n\}_{n \ge 0}$ and $\{ \alpha n \}_{n \ge
0}$ have the joint repetition property. This implies that
$\{\omega_n\}_{n \ge 0}$ and $\{\omega_n + \alpha n\}_{n \ge 0}$
have the joint repetition property.

Clearly, the addition of $\beta$ does not affect the distances in
question so that the asserted joint repetition property for
$\{\omega_n\}_{n \ge 0}$ and $\{\omega_n + \alpha n + \beta\}_{n \ge
0}$ follows immediately.
\end{proof}

Note again that there is an immediate extension of the first part
of the lemma to sequences in a general metric space that have the
repetition property. Any such sequence has the joint repetition
property with the sequence $\{ \alpha n \}_{n \ge 0}$ in $\T$ for
any $\alpha \in \T$. The result in this generality will be used in
the proof of Theorem~\ref{t.ap} below.

The final preservation property we wish to address here is the
multiplication of a given sequence by a rational number. The
following lemma exhibits a class of sequences for which the
repetition property is preserved under such an operation.

\begin{lemma}\label{l.recurrence}
Suppose $\{\omega_n\}_{n \ge 0}$ is a sequence in $\R$ whose
projection to $\T$ has the repetition property and which satisfies a
linear recurrence relation of the following form: there exist
integers $k \ge 1$ and $a_1, \ldots, a_k$ with $|a_k| = 1$ such that
$$
\omega_n + a_1 \omega_{n-1} + a_2 \omega_{n-2} + \cdots + a_k
\omega_{n-k} = 0.
$$
Then, for every $r \in \Q$, the projection of the sequence $\{r
\omega_n\}_{n \ge 0}$ to $\T$ has the repetition property, jointly
with the original sequence.
\end{lemma}

\begin{proof}
By Lemma~\ref{l.basic}.(c), it suffices to consider the case $r =
m^{-1}$, where $m$ is an integer $\ge 2$. By the recurrence
relation, any $k$ consecutive values determine the entire sequence.
This is of course also true for the sequence $\{m^{-1} \omega_n\}_{n
\ge 0}$ in $\R$.

By assumption, $\{\omega_n \!\! \mod 1 \}_{n \ge 0}$ has the
repetition property in $\T = \R / \Z$. Consequently, $\{m^{-1}
\omega_n \!\! \mod m^{-1} \}_{n \ge 0}$ has the repetition property
in $\R / (m^{-1} \Z)$ with, in fact, the same $q(\varepsilon,r)$ as
the original sequence. To see that $\{m^{-1} \omega_n \!\! \mod 1
\}_{n \ge 0}$ has the repetition property, notice that when passing
from $\R / (m^{-1} \Z)$ to $\R / \Z$, we can in principle draw from
$m$ choices for each $n$. After fixing $k$ values, the rest of the
sequence is effectively determined.\footnote{It is determined only
up to integers, which however are irrelevant for the repetition
property in $\T$.}

Now consider almost-repetitions of $m^{-1} \omega_n \!\! \mod
m^{-1}$ in $\R / (m^{-1} \Z)$ with $r$ sufficiently large (larger
than $m^k$) and apply the pigeonhole principle to see that $m^{-1}
\omega_n \!\! \mod 1$ in $\R / \Z$ must have almost-repetitions as
well. More precisely, find a repeated cell of size $m^{-1}$ first
and then use the recursion relation to see that the resulting
pattern must repeat in both directions.

In this way, we can map every triple
$(\varepsilon,r,q(\varepsilon,r))$ describing an almost-repetition
for the original sequence $\{\omega_n \!\! \mod 1 \}_{n \ge 0}$ to
a new triple for the derived sequence $\{m^{-1} \omega_n \!\! \mod
1 \}_{n \ge 0}$. With these derived almost-repetitions, one may
check that $\{m^{-1} \omega_n \!\! \mod 1 \}_{n \ge 0}$ has the
repetition property in $\T$, jointly with the projection of
$\{\omega_n\}_{n \ge 0}$ to $\T$.
\end{proof}

\subsection{Characterization of Polynomials With the
Repetition Property}

Consider the sequence $\omega_n = p(n) \!\! \mod 1$ in $\T$, where
$p(n) = \sum_{k=0}^d a_k n^k$ is a polynomial. If this sequence
has the repetition property, we say that $p$ has the repetition
property. Let us denote the terms of $p$ by $p_k$, that is,
$p_k(n) = a_k n^k$.

\begin{theorem}\label{thm.ploychar}
A polynomial $p$ has the repetition property if and only if its
terms $p_k$, $0 \le k \le d$ have the joint repetition property.
\end{theorem}

\begin{proof}
That $p$ has the repetition property if $p_k$, $0 \le k \le d$ have
the joint repetition property is easy to verify.

Assume now that $p$ has the repetition property. By
Lemma~\ref{l.basic}.(b), the polynomials $p^{(l)}$ given by
$p^{(l)}(n) = p (n l)$, $0 \le l \le d$ have the joint repetition
property. Moreover, by taking a finite number of discrete
derivatives of $p$, we eventually obtain the zero polynomial. In
this way, we find a linear recursion relation between $p(n) , \ldots
, p(n-d-1)$ with constant coefficients whose absolute values are
given by the entries in the Pascal triangle. Consequently, $p(n)$
obeys the assumptions of Lemma~\ref{l.recurrence}.

Consider the matrix $A = (a_{i,j})_{0 \le i,j \le d}$ given by
$a_{i,j} = i^j$. Then
$$
\begin{pmatrix} p(n \cdot 0) \\ p(n \cdot 1) \\ p( n \cdot 2) \\
\vdots \\ p(n \cdot d) \end{pmatrix} = A \begin{pmatrix} a_0 \\ a_1 n \\ a_2 n^2 \\
\vdots \\ a_d n^d \end{pmatrix} .
$$
Notice that $A$ a Vandermonde matrix and hence its invertibility
follows from the well-known formula for the determinant of a
Vandermonde matrix. Of course, $A^{-1}$ has only rational entries.

Since
$$
\begin{pmatrix} a_0 \\ a_1 n \\ a_2 n^2 \\ \vdots \\ a_d n^d
\end{pmatrix} = A^{-1} \begin{pmatrix} p(n \cdot 0) \\ p(n \cdot 1)
\\ p( n \cdot 2) \\ \vdots \\ p(n \cdot d) \end{pmatrix},
$$
it therefore follows from Lemma~\ref{l.recurrence} that $p_k$, $0
\le k \le d$ have the joint repetition property.
\end{proof}

Theorem~\ref{thm.ploychar} suggests investigating the repetition
property for homogeneous polynomials. This problem is addressed in
the following theorem.

\begin{theorem}\label{t.hompoly}
A homogeneous polynomial $p_k$ of the form $p_k(n) = a_k n^k$ has
the repetition property if and only if $\liminf_{q \to \infty}
q^{k-1} \langle a_k q \rangle = 0$.
\end{theorem}

\begin{proof}
Write $\omega_n = p_k(n) \!\! \mod 1$. We have
$$
\mathrm{dist}(\omega_{n + q},\omega_n) = \left\langle \sum_{j =
0}^{k - 1} a_k
\begin{pmatrix} k \\ j \end{pmatrix} n^j q^{k-j} \right\rangle
$$
This shows immediately that the repetition property follows from the
existence of $q_m \to \infty$ with $\lim_{m \to \infty} q_m^{k-1}
\langle a_k q_m \rangle = 0$.

Conversely, assume that $p_k$ has the repetition property. By
Lemma~\ref{l.basic}.(a), $\tilde p_k$ given by
$$
\tilde p_k(n) = p_k(n+1) = a_k \sum_{j = 0}^k
\begin{pmatrix} k \\ j \end{pmatrix} n^j
$$
has the repetition property as well. Thus, by
Theorem~\ref{thm.ploychar} and Lemma~\ref{l.recurrence}, we find
that $n \mapsto a_k n^j$, $0 \le j \le k$ have the joint
repetition property. For $j = 0$ and $j = 1$, this statement is
obvious. Let us explore what information can be gleaned from the
case $j = 2$: for every $\varepsilon > 0$, we have $\left\langle 2
a_k q_l n + a_k q_l^2 \right\rangle < \varepsilon$ for some $q_l
\to \infty$ and $0 \le n \le q_l$. Evaluating this for $n = 0$, we
find that $\langle a_k q_l^2 \rangle < \varepsilon$. Now vary $n$.
Each time we increase $n$, we shift in the same direction by
$\langle 2 a_k q_l \rangle$. If $\varepsilon > 0$ is sufficiently
small, it follows that we cannot go around the circle completely
and hence we have $\langle 2n a_k q_l \rangle = n \langle 2 a_k
q_l \rangle$ for every $0 \le n \le q_l$. Taking $\varepsilon$ to
zero, this implies $\liminf_{q \to \infty} q \langle a_k q \rangle
= 0$. Using this, we can consider the case $j = 3$ and argue in a
similar way to find that $\liminf_{q \to \infty} q^2 \langle a_k q
\rangle = 0$. Carrying on inductively, we arrive at the desired
statement, $\liminf_{q \to \infty} q^{k-1} \langle a_k q \rangle =
0$.
\end{proof}

For the remainder of this section, we will be concerned with a
special class of polynomials; namely, $p(n) = \alpha n^3 + \beta
n^2$. The study of the repetition property for these polynomials
is somewhat more involved as they present a borderline case. We
have the following theorem:

\begin{theorem}\label{t.border}
There are a residual subset $P_0$ of $\T$ and a subset $G$ of $\T$
of full Lebesgue measure such that for $\alpha \in P_0$ and $\beta
\in G$, the polynomial $p(n) = \alpha n^3 + \beta n^2$ has the
repetition property.
\end{theorem}

\begin{proof}
Set $X = \T^2$. Lemma~\ref{l.gdelta} implies that
$$
H = \{ (\alpha , \beta) \in \T^2 : p(n) = \alpha n^3 + \beta n^2
\text{ does not have the repetition property} \, \}
$$
is an $F_\sigma$ subset of $\T^2$.

Next we wish to apply Lemma~\ref{l.chebosh}. We set $P = K = \T$,
so that $W = P \times K = \T^2$, and let $L$ be any $G_\delta$
subset of $\T$ of Lebesgue measure zero that contains the set of
badly approximable numbers (those numbers that have bounded
partial quotients).

Writing
\begin{align*}
K(\alpha) & = \{ \beta \in \T : (\alpha,\beta) \in H \} \\
& = \{ \beta \in \T : p(n) = \alpha n^3 + \beta n^2 \text{ does
not have the repetition property} \, \}
\end{align*}
for $\alpha \in \T$ and
$$
P_0 = \{ \alpha \in \T : K(\alpha) \subset L \},
$$
and observing that $P_0$ contains all rational numbers in $\T$, it
follows from Lemma~\ref{l.chebosh} that $P_0$ is a residual subset
of $\T$. In other words, if we denote $G = \T \setminus L$, then
$p(n) = \alpha n^3 + \beta n^2$ has the repetition property for
every $\alpha \in P_0$ and $\beta \in G$.
\end{proof}

\section{Proof of Theorem~\ref{t.skewshift}}\label{s.skew}

Denote the $j$-th component of $T^n \omega$ by $t(n,j,\omega)$.
Observe that $n \mapsto T^n \omega$ has the repetition property if
and only if $n \mapsto t(n,j,\omega)$, $1 \le j \le k$ have the
joint repetition property.

Clearly, we have
\begin{equation}\label{tn1}
t(n,1,\omega) = n \alpha + \omega_1.
\end{equation}
Noting that $t(n,2,\omega) = t(n-1,1,\omega) + t(n-1,2,\omega)$,
we find that $t(n,2,\omega) = t(n-1,1,\omega) + t(n-2,1,\omega) +
\cdots + t(0,1,\omega) + t(0,2,\omega)$ and hence
\begin{equation}\label{tn2}
t(n,2,\omega) = \frac{n(n-1)}{2} \alpha + n \omega_1 + \omega_2.
\end{equation}
Continuing in this fashion, we see that $t(n,j,\omega)$ is a
polynomial in $n$ of degree $j$ and its leading coefficient is a
rational multiple of $\alpha$, while the other coefficients are
rational linear combinations of $\alpha,\omega_1,\ldots,\omega_k$,
but not pure multiples of $\alpha$.
\\[1mm]
(a) Assume that $T$ has (TRP). In particular, for some $\omega \in
\T^k$, the polynomial $n \mapsto t(n,k,\omega)$ has the repetition
property. Thus, by Theorem~\ref{thm.ploychar} and
Lemma~\ref{l.recurrence}, the polynomial $n \mapsto \alpha n^k$
has the repetition property, which is equivalent to $\liminf_{q
\to \infty} q^{k-1} \langle q \alpha \rangle = 0$ by
Theorem~\ref{t.hompoly}. Conversely, if the polynomial $n \mapsto
\alpha n^k$ has the repetition property, this is true jointly for
all polynomials $n \mapsto \alpha n^j$, $0 \le j \le k$ (by the
argument in the proof of Theorem~\ref{t.hompoly}). It follows that
$0 \in PRP(T)$ and hence $T$ has (TRP).
\\[1mm]
(b) This was shown in \cite{bd1} and is stated here for the sake
of completeness.
\\[1mm]
(c) If $k \ge 3$, then $n \mapsto t(n,3,\omega)$ is a polynomial
of degree $3$ and the coefficient of $n^2$ is a rational linear
combination of $\alpha,\omega_1,\ldots,\omega_k$, but not a pure
multiple of $\alpha$. In particular, for some $\omega \in \T^k$,
this coefficient is badly approximable and hence, by
Theorems~\ref{thm.ploychar} and \ref{t.hompoly}, this $\omega$
does not belong to $PRP(T)$. In other words, $T$ does not have
(GRP).
\\[1mm]
(d) This follows immediately from part (a).
\\[1mm]
(e) We begin with $t(n,3,\omega)$, discuss the repetition property
for this sequence in the circle, and then proceed to include the
other components of $T^n \omega$ in our consideration. As we saw
in Lemma~\ref{l.addition}, we can restrict our attention to the
terms in $t(n,3,\omega)$ of degree $2$ and $3$. Computing these
terms, we find that
$$
t(n,3,\omega) = \frac{\alpha}{6} n^3 + \left( -\frac{\alpha}{2} +
\frac{\omega_1}{2} \right) n^2 + \text{ lower order terms}.
$$
Applying Theorem~\ref{t.border}, we find that $t(n,3,\omega)$ has
the repetition property for $\alpha$'s from a certain residual
subset of $\T$ and $\omega_1$'s from a certain full measure subset
of $\T$. Recalling the specific form \eqref{tn2} of
$t(n,2,\omega)$, we can use Theorems~\ref{thm.ploychar} and
\ref{t.hompoly} to see that $t(n,2,\omega)$ has the repetition
property, jointly with $t(n,3,\omega)$, for $\alpha$ and
$\omega_1$ as above and every $\omega_2$. Since $t(n,1,\omega)$
has the form \eqref{tn1}, adding on the joint repetition property
for $t(n,1,\omega)$ is then easy by the argument from the proof of
Lemma~\ref{l.addition}. Putting everything together, we find that
for $\alpha$'s from a certain residual subset of $\T$,
$\omega_1$'s from a certain full measure subset of $\T$, and every
$(\omega_2,\omega_3) \in \T^2$, $t(n,j,\omega)$, $1 \le j \le 3$
have the joint repetition property. Consequently, $T$ has (MRP)
for such $\alpha$'s.
\\[1mm]
(f) If $k \ge 4$, then $n \mapsto t(n,4,\omega)$ is a polynomial
of degree $4$ and the coefficient of $n^3$ is a rational linear
combination of $\alpha,\omega_1,\ldots,\omega_k$, but not a pure
multiple of $\alpha$. Consider an $\omega_i$ with non-zero
coefficient in this linear combination. Fixing all the other
entries of $\omega$, we see that the coefficient satisfies the
necessary condition for the repetition property in this component
only for $\omega_i$ from a set of zero Lebesgue measure. By
Fubini, it follows that $T$ does not have (MRP).\hfill\qedsymbol

\begin{appendix}

\section{Universal Joint Repetition Property}

The notion of joint repetition property was useful in our study
above. Pushing a bit further, it is natural to introduce the
following class of sequences:

\begin{definition}
A sequence $\{ \omega_n \}_{n \ge 0}$ in a metric space $\Omega$ has
the universal joint repetition property if for every sequence
$\{\tilde \omega_n\}_{n \ge 0}$ in some metric space $\tilde \Omega$
that has the repetition property, $\{\omega_n\}_{n \ge 0}$ and
$\{\tilde \omega_n\}_{n \ge 0}$ have the joint repetition property.
\end{definition}

We consider it an interesting problem to find a characterization
of the sequences that have the universal joint repetition
property. Theorem~\ref{t.ap} below presents a first step in this
direction.

Recall that $a : \Z \to \C$ is called almost periodic if its
translates form a relatively compact subset of $\ell^\infty(\Z)$.
It is a fundamental result (see, e.g., \cite{b1}) that every
almost periodic sequence can be approximated uniformly by finite
linear combinations of exponentials $e^{2 \pi i \alpha n}$. The
following theorem is an extension of Lemma~\ref{l.addition} in
Section~\ref{s.poly}.

\begin{theorem}\label{t.ap}
Every almost periodic sequence has the universal joint repetition
property.
\end{theorem}

\begin{proof}
Let $\{ \omega_n\}_{n \ge 0}$ be a sequence in a metric space
$\Omega$ that has the repetition property. By
Lemma~\ref{l.addition} (see the remark after the proof of this
lemma), $\{ \omega_n\}_{n \ge 0}$ and $\{ \alpha n \}_{n \ge 0}$
in $\T$ have the joint repetition property for every $\alpha \in
\T$. Thus, by Lemmas~\ref{l.joint1} and \ref{l.joint2}, $\{
\omega_n\}_{n \ge 0}$ and $\{ c e^{2 \pi i \alpha n} \}_{n \ge 0}$
in $\C$ have the joint repetition property for every $\alpha \in
\T$ and $c \in \C$.

Inductively, it follows that $\{ \omega_n\}_{n \ge 0}$ and any
finite linear combination of exponentials have the joint
repetition property. Since the latter lie densely in the almost
periodic sequences with respect to uniform convergence, the
theorem now follows from Lemma~\ref{l.joint3}.
\end{proof}

\end{appendix}

\end{document}